\documentclass{article}

\usepackage{amsmath}
\usepackage{amssymb}
\usepackage{amsthm}
\usepackage{authblk}
\usepackage{url}

\usepackage[style=numeric,sorting=none]{biblatex}
\newcommand{\sd}{\mathrm{d}}

\DeclareMathOperator{\sgn}{sgn}

\title{Fourier transform of composed functions}
\author[1]{David Venhoek}%\footnote{Work for this paper was done at the Theoretical High Energy Physics department at Radboud University.}}
\affil[1]{Theoretical High Energy Physics, Radboud University
Nijmegen, Heyendaalseweg 135, 6525~AJ Nijmegen, The Netherlands}
\date{\today}

\allowdisplaybreaks

\theoremstyle{plain}
\newtheorem{lemma}{Lemma}
\newtheorem{theorem}{Theorem}

\theoremstyle{definition}
\newtheorem{definition}{Definition}
\newtheorem{notation}{Notation}

\begin{document}

\maketitle

\begin{abstract}
We prove an explicit formula for the Fourier transform of $f(u(t))$, given the Fourier transform of $f(t)$, assuming $f\in L^2(-\infty,\infty)$ and $u$ sufficiently well behaved. We illustrate its usefulness by calculating the Fourier transform of $(a^2 + \sinh(bt)^2)^{-1}$.
\end{abstract}

The interaction between the Fourier transform and function composition has applications in Computer Graphics and sampling theory\cite{Bergner}. Furthermore, a similar pattern of function composition occurs in Physics, particularly in the study of Unruh detectors\cite{Unruh, DeWitt, beenakker2023accelerating}.

Bergner et al.~\cite{Bergner} give a short derivation for an explicit formula for calculating the Fourier transform of the composed function in terms of the Fourier transform of the outer function. However, their derivation is not mathematically sound as it involves a change of integration order that cannot be justified by Fubini's theorem. This paper provides a proper proof showing the formula from Bergner et al. is valid for $f(u(t))$ when $f\in L^2(-\infty,\infty)$ and $u$ sufficiently well behaved.

\section{Conventions and notation}

Throughout the following we will encounter several improper integrals. For brevity of notation, we will follow the convention in Titchmarsh~\cite{Titchmarsh} for these:
\begin{notation}[Improper integrals]
For one-sided improper integrals, we write
\begin{align*}
\int\limits_a^{\rightarrow\infty} f(t)\sd t &\equiv \lim_{T\rightarrow\infty} \int\limits_a^T f(t)\sd t,
\end{align*}
which we generalize for integrals with both endpoints improper to
\begin{align*}
\int\limits_{\rightarrow -\infty}^{\rightarrow \infty} f(t)\sd t &\equiv \lim_{S,T\rightarrow\infty} \int\limits_{-S}^T f(t)\sd t.
\end{align*}
\end{notation}
Note that if $f(t)$ is Lebesgue integrable on $\mathbb{R}$, the above limit converges to the Lebesgue integral of $f$.

\begin{notation}[Fourier transforms]
For a function $f\in L^2(-\infty,\infty)$ we denote its Fourier transform with $\hat{f}(k) = \int\limits_{\rightarrow-\infty}^{\rightarrow\infty} e^{ikt}f(t)\sd t$. In some cases, we will define $\hat{f}$ first, in which case $f$ will be the corresponding function whose Fourier transform is $\hat{f}$. By Plancherel's theorem this exists so long as $\hat{f}\in L^2(-\infty,\infty)$.
\end{notation}

Throughout this text, when considering function spaces, we will always endow these spaces with the $L^2$ norm. For norms and limits in these spaces we define the following shorthand notation:
\begin{notation}[$L^2$ convergence]
Let $\{f_n\}$ be a sequence of functions in $L^2(-\infty,\infty)$. We write $||f_n||$ for the $L^2$ norm of $f_n$. We write $f_n\rightarrow f$ if the sequence of functions $f_n$ converges to $f$ in the $L^2$ norm.
\end{notation}

In the context of Fourier theory, Schwartz space occupies a special role in that the Fourier transform integral is integrable for all frequencies. Furthermore, the Fourier transform is a bijection from Schwartz space to itself.
\begin{definition}[Schwartz space]
Let $\mathcal{S}$ denote Schwartz space. That is, the space of smooth functions $f$ with $\sup\limits_{x\in\mathbb{R}} |x^\alpha f^{(\beta)}(x)| < \infty$ for all natural numbers $\alpha$, $\beta$.
\end{definition}
Combined with Plancherel's theorem, which states that the Fourier transform preserves the $L^2$ norm, and the fact that Schwartz space is dense as a subset of $L^2$, this extends the Fourier transform to a unitary map on $L^2$.\cite{stein2011fourier, stein2009real}

Finally, we will have need to consider functions that show various types of asymptotic behaviour on the extremes of their domain:
\begin{definition}[Proper functions]
A function $u:\mathbb{R}\rightarrow\mathbb{R}$ is \emph{proper} iff $\left|u(t)\right|\rightarrow \infty$ when $t\rightarrow\pm\infty$.
\end{definition}
\begin{definition}[Eventual monotonicity]
A function $u:\mathbb{R}\rightarrow\mathbb{R}$ is \emph{eventually monotone} iff there exists $S, T \in\mathbb{R}$ such that the restrictions $u_S: (-\infty,S)\rightarrow \mathbb{R}$ and $u_T: (T,\infty)\rightarrow\mathbb{R}$ are monotone.
\end{definition}
Note that we do not require monotone functions to be injective.

For example, the function $u(t) = |t|$ is both proper and eventually monotone, and so is $u(t) = t^3 + sin(t)$.

\section{Transfer functions}

Let us start by a study of the transfer function that we will find for the Fourier transformation of function composition. For a function $u\in C^1\left(\mathbb{R}\right)$ with $u'$ proper and eventually monotone, we define the Fresnel-type integral
\begin{align}
H_u(k,l) &\equiv \int\limits_{\rightarrow-\infty}^{\rightarrow\infty} e^{i(kt-lu(t))}\sd t.
\end{align}

\begin{lemma}\label{transfer}
If $u\in C^1\left(\mathbb{R}\right)$ with $u'$ proper and eventually monotone, then for any $\epsilon>0$ the limit defining $H_u(k,l)$ converges uniformly in $l$ on the set $\mathbb{R}\setminus\left(-\epsilon,\epsilon\right)$.
\end{lemma}

\begin{proof}
The definition of $H_u(k,l)$ can be split
\begin{align*}
H_u(k,l) &= \int\limits_{\rightarrow-\infty}^0 e^{i(kt-lu(t))}\sd t + \int\limits_0^{\rightarrow\infty} e^{i(kt-lu(t))}\sd t,
\end{align*}
and as a result of the symmetry of this expression it is sufficient to show uniform convergence for
\begin{align*}
\int\limits_0^{\rightarrow\infty} e^{i(kt-lu(t))}\sd t.
\end{align*}
Furthermore, without loss of generality we can assume $l>0$ and $u'(t)\rightarrow +\infty$ as $t\rightarrow +\infty$.

Now let $M$ be sufficiently large so that $u'(t)$ is monotone on $[M,\infty)$ and $\left|\epsilon u'(M)\right| > \left|k\right|$. Then $v(t) = kt-lu(t)$ is a continuous injection on $[M,\infty)$ and hence
\begin{align*}
\left|\int\limits_0^{\rightarrow\infty} e^{i(kt-lu(t))}\sd t - \int\limits_0^M e^{i(kt-lu(t))}\sd t\right| &= \left|\int\limits_M^{\rightarrow\infty} e^{i(kt-lu(t))}\sd t\right|\\
&= \left|\int\limits_{kM-lu(M)}^{\rightarrow-\infty} \frac{e^{iv}}{k-lu'(t(v))}\sd v\right|\\
&\le \frac{4\pi}{\left|k-lu'(M)\right|}\\
&\le \frac{4\pi}{\left|\epsilon u'(M)\right| - \left|k\right|},
\end{align*}
where the second to last step follows from monotonicity of $u'(t)$ and hence of $k-lu'(t)$. Since $u'(M)\rightarrow \infty$ as $M\rightarrow \infty$, and with the last expression independent of $l$, this implies that the limit converges uniformly in $l$. By the previous observations this implies that the defining limit of $H_u(k,l)$ exists and converges uniformly for $l \in \mathbb{R}\setminus\left(\epsilon,-\epsilon\right)$.
\end{proof}

\section{Fourier transform of compositions}
We can now follow the standard approach for proving facts about the Fourier transform on $L^2$. We start by proving our desired result on a subspace of the Schwartz functions that is (under the $L^2$ norm) dense in all of Schwartz space, and then use limiting procedures and Plancherel's theorem to extend it to the entirety of $L^2(-\infty,\infty)$.

\begin{lemma}\label{subcase}
Let $f\in\mathcal{S}$, with $\hat{f}(\omega) = 0$ when $|\omega| < \mu$ for some $\mu > 0$, and let $u \in C^1$ with $u'$ proper and eventually monotone. Then
\begin{align}
\int\limits_{-\infty}^{\infty} e^{ikt}f(u(t))\sd t = \frac{1}{2\pi}\int\limits_{-\infty}^{\infty}\hat{f}(l)H_u(k,l)\sd l.
\end{align}
\end{lemma}

\begin{proof}
Calculating, we find
\begin{align*}
\int\limits_{-\infty}^{\infty} e^{ikt}f(u(t))\sd t &= \lim_{S,T\rightarrow\infty} \int\limits_{-S}^T e^{ikt} f(u(t))\sd t\\
&=\lim_{S,T\rightarrow\infty} \int\limits_{-S}^T e^{ikt}\frac{1}{2\pi}\int\limits_{-\infty}^\infty\hat{f}(l)e^{-ilu(t)}\sd l\sd t\\
&=\lim_{S,T\rightarrow\infty} \frac{1}{2\pi} \int\limits_{-\infty}^\infty\hat{f}(l)\int\limits_{-S}^T e^{i(kt-lu(t))}\sd t\sd l\\
&= \frac{1}{2\pi}\int\limits_{-\infty}^\infty\hat{f}(l)H_u(k,l)\sd l.
\end{align*}
The second step above is justified by Fubini's theorem. To see the final step is justified note that because $f\in\mathcal{S}$, also $\hat{f}\in\mathcal{S}$. Hence, by definition of Schwartz space $\hat{f}$ is bounded. This, combined with the fact that $\hat{f}$ is supported in $\mathbb{R}\setminus(-\mu,\mu)$ implies Lemma~\ref{transfer} applies and shows that the integrand converges uniformly. This completes the proof. 
\end{proof}

In applying the standard limiting procedure we will have a sequence of functions $f_n$ converging to $f$. But for this to work well with the above result, we then also need control over the convergence of $f_n\circ u$:

\begin{lemma}\label{convergence}
Let $f_n \in L^2(-\infty,\infty)$ with $f_n \rightarrow f$, and and $u:\mathbb{R}\rightarrow \mathbb{R}$ a bijection with $u\in C^1$ and $|u'(t)| > C$ for all $t\in\mathbb{R}$. Then $f_n\circ u \rightarrow f\circ u$.
\end{lemma}

\begin{proof}
Calculating, we find
\begin{align*}
||f_n\circ u - f\circ u||^2 &= \int\limits_{-\infty}^\infty |f_n(u(t)) - f(u(t))|^2 \sd t\\
&= \int\limits_{-\infty}^\infty |f_n(u) - f(u)|^2 \frac{1}{|u'(t(u))|} \sd u\\
&\le \frac{1}{C} ||f_n - f||^2.
\end{align*}
Since $f_n\rightarrow f$, this last expression can be made arbitrarily small by increasing $n$. Hence $f_n\circ u\rightarrow f\circ u$.
\end{proof}

With these ingredients, we are ready for our main result.

\begin{theorem}\label{mainresult}
Let $f\in L^2(-\infty,\infty)$, and $u:\mathbb{R}\rightarrow\mathbb{R}$ a bijection with $u\in C^1$ and $|u'(t)| > C$ for all $t\in\mathbb{R}$. Then
\begin{align}
\int\limits_{\rightarrow-\infty}^{\rightarrow\infty} e^{ikt}f(u(t))\sd t = \frac{1}{2\pi}\int\limits_{\rightarrow-\infty}^{\rightarrow\infty}\hat{f}(l)H_u(k,l)\sd l
\end{align}
for almost every $k$.
\end{theorem}

\begin{proof}
Since Schwartz functions are dense in $L^2$, so are Schwartz functions that are zero in a region around zero. Hence we can find a sequence of functions $\hat{f}_n$ in $\mathcal{S}$ such that $\hat{f}_n\rightarrow \hat{f}$ and there exist $\epsilon_n > 0$ with $\hat{f}_n(k) = 0$ when $|k| < \epsilon_n$. Since the Fourier transform is a bijection on $\mathcal{S}$, we have $f_n \in \mathcal{S}$, and furthermore by Plancherel's theorem $f_n \rightarrow f$.

Now, define
\begin{align*}
g_n(k) = \frac{1}{2\pi}\int\limits_{-\infty}^{\infty}\hat{f}_n(l)H_u(k,l)\sd l,
\end{align*}
then by Lemma~\ref{subcase} we have $\widehat{f_n\circ u} = g_n$. Furthermore, by Lemma~\ref{convergence}, $f_n\circ u \rightarrow f\circ u$.

Applying Plancherel one final time, this implies $\widehat{f_n\circ u} \rightarrow \widehat{f\circ u}$, and hence $g_n\rightarrow \widehat{f\circ u}$, proving the theorem.
\end{proof}

\section{Application: Fourier transform of $\frac{1}{a^2+\sinh(bt)^2}$}

The above result opens a new avenue for calculating the Fourier transform of functions. We will illustrate this by using it to establish
\begin{align}\label{example_target}
\hat{g}(k) &= \frac{2\pi\cosh\left(\frac{k\pi}{2b}\right)}{ab\sinh\left(\frac{k\pi}{b}\right)} \frac{\sinh\left(\frac{k}{b}\arccos\left(a\right)\right)}{\sqrt{1-a^2}}
\end{align}
as the Fourier transform of
\begin{align}
g(t) &= \frac{1}{a^2+\sinh(bt)^2},
\end{align}
assuming both $a$ and $b$ positive, taking the principal branches for the values of the inverse cosine and square roots.

Note that, to the best of the knowledge of the authors, this result is not present in the standard sources for integrals such as \cite{GradRyzh} or \cite{DLMF}, although there may be less well-known sources that contain it. Furthermore, the approach shown here is likely not the only viable way of calculating the desired result.

We will start by defining
\begin{align*}
f(t) &= \frac{1}{a^2+t^2},\\
u(t) &= \sinh(bt),
\end{align*}
and noting that this implies that $g(t) = f(u(t))$.

A quick lookup in for example~\cite{GradRyzh} shows that the Fourier transform of $f$ is given by
\begin{align*}
\hat{f}(k) &= \frac{\pi e^{-a\left|k\right|}}{a}.
\end{align*}
The function $H_u$ can be calculated using formula 10.32.7 in~\cite{DLMF} and some straightforward manipulations to be
\begin{align*}
H_u(k,l) &= \frac{2}{b}e^{\sgn(l)\frac{k\pi}{2b}}K_{i\frac{k}{b}}\left(\left|l\right|\right),
\end{align*}
where $\sgn(l)$ denotes the sign of $l$ (that is, $+1$ if $l$ is positive, $-1$ if $l$ is negative and $0$ otherwise).
Combining this, Theorem~\ref{mainresult} then shows that
\begin{align*}
\hat{g}(k) &= \frac{2\cosh\left(\frac{k\pi}{2b}\right)}{ab}\int_0^{\infty} e^{-al}K_{i\frac{k}{b}}\left(l\right)\sd l.
\end{align*}

The integral on the right hand side can be found in~\cite{GradRyzh}, equation 6.611.3. Applying this, a straightforward calculation gives the result in Equation~\ref{example_target}.

\section{Acknowledgements}

The author would like to thank Dr. M\"uger for his input on early drafts of this paper.

\printbibliography
\end{document}